
\font \tenbfne                = cmbx10
\font \eightbf                = cmbx8
\font \eightbfne              = cmbx10 at 8pt
\font \eighti                 = cmmi8
\font \eightit                = cmti8
\font \eightrm                = cmr8
\font \eightsans              = cmss10 at 8pt
\font \eightsl                = cmsl8
\font \eightsy                = cmsy8
\font \eighttt                = cmtt8
\font \fivesans               = cmss10 at 5pt
\font \kk                     = cmmi6

\font \markfont               = cmr10 scaled\magstep1
\font \subhfont               = cmr10 scaled\magstep4
\font \sevensans              = cmss10 at 7pt
\font \sixbf                  = cmbx6
\font \sixbfne                = cmbx10 at 6pt
\font \sevenbfne              = cmbx10 at 7pt
\font \fivebfne               = cmbx10 at 5pt
\font \sixi                   = cmmi6
\font \sixrm                  = cmr6
\font \authfont               = cmr17 scaled\magstep2
\font \headfont               = cmbx12 scaled\magstep4
\font \sixsans                = cmss10 at 6pt
\font \sixsy                  = cmsy6
\font \smallescriptfont       = cmr5 at 7pt
\font \smallescriptscriptfont = cmr5
\font \smalletextfont         = cmr5 at 10pt
\font \tafontt                = cmbx10 scaled\magstep2
\font \tafonts                = cmbx7  scaled\magstep2
\font \tafontss               = cmbx5  scaled\magstep2
\font \tamt                   = cmmib10 scaled\magstep2
\font \tams                   = cmmib10
\font \tamss                  = cmmib10 scaled 700
\font \tast                   = cmsy10 scaled\magstep2
\font \tass                   = cmsy7  scaled\magstep2
\font \tasss                  = cmsy5  scaled\magstep2
\font \tasyt                  = cmex10 scaled\magstep2
\font \tasys                  = cmex10 scaled\magstep1
\font \tbfontt                = cmbx10 scaled\magstep1
\font \tbfonts                = cmbx7  scaled\magstep1
\font \tbfontss               = cmbx5  scaled\magstep1
\font \tbmt                   = cmmib10 scaled\magstep1
\font \tbms                   = cmmib10 scaled 833
\font \tbmss                  = cmmib10 scaled 600
\font \tbst                   = cmsy10 scaled\magstep1
\font \tbss                   = cmsy7  scaled\magstep1
\font \tbsss                  = cmsy5  scaled\magstep1
\font \tensans                = cmss10
\magnification=\magstep1
\vsize=19.4cm
\hsize=14.4truecm
\hfuzz=2pt
\tolerance=10000
\abovedisplayskip=3 mm plus6pt minus 4pt
\belowdisplayskip=3 mm plus6pt minus 4pt
\abovedisplayshortskip=0mm plus6pt minus 2pt
\belowdisplayshortskip=2 mm plus4pt minus 4pt
\predisplaypenalty=0
\clubpenalty=10000
\widowpenalty=10000
\frenchspacing
\parindent=1.5em
\newdimen\oldparindent\oldparindent\parindent
\def\bf{\fam\bffam\tenbfne}
\textfont\bffam=\tenbfne \scriptfont\bffam=\sevenbfne
\scriptscriptfont\bffam=\fivebfne
\skewchar\eighti='177 \skewchar\sixi='177
\skewchar\eightsy='60 \skewchar\sixsy='60
\hyphenchar\eighttt=-1
\def\newline{\hfil\break}%
\newtoks\eq\newtoks\eqn
\catcode`@=11
\def\folio{\ifnum\pageno<\z@
\uppercase\expandafter{\romannumeral-\pageno}%
\else\number\pageno \fi}
\def\eqalign#1{\null\vcenter{\openup\jot\m@th
  \ialign{\strut\hfil$\displaystyle{##}$&$\displaystyle{{}##}$\hfil
      \crcr#1\crcr}}}
\def\displaylines#1{{}$\displ@y
\hbox{\vbox{\halign{$\@lign\hfil\displaystyle##\hfil$\crcr
    #1\crcr}}}${}}
\def\eqalignno#1{{}$\displ@y
  \hbox{\vbox{\halign to\hsize{\hfil$\@lign\displaystyle{##}$\tabskip\z@skip
    &$\@lign\displaystyle{{}##}$\hfil\tabskip\centering
    &\llap{$\@lign##$}\tabskip\z@skip\crcr
    #1\crcr}}}${}}
\def\leqalignno#1{{}$\displ@y
\hbox{\vbox{\halign
to\hsize{\qquad\hfil$\@lign\displaystyle{##}$\tabskip\z@skip
    &$\@lign\displaystyle{{}##}$\hfil\tabskip\centering
    &\kern-\hsize\rlap{$\@lign##$}\tabskip\hsize\crcr
    #1\crcr}}}${}}
\def\generaldisplay{%
\ifeqno
       \ifleqno\leftline{$\displaystyle\the\eqn\quad\the\eq$}%
       \else\line{\kern\parindent$\displaystyle\the\eq\hfill\the\eqn$}\fi
\else
       \leftline{\kern\parindent$\displaystyle\the\eq$}%
\fi
\global\eq={}\global\eqn={}}%
\newif\ifeqno\newif\ifleqno \everydisplay{\displaysetup}
\def\displaysetup#1$${\displaytest#1\eqno\eqno\displaytest}
\def\displaytest#1\eqno#2\eqno#3\displaytest{%
\if!#3!\ldisplaytest#1\leqno\leqno\ldisplaytest
\else\eqnotrue\leqnofalse\eqn={#2}\eq={#1}\fi
\generaldisplay$$}
\def\ldisplaytest#1\leqno#2\leqno#3\ldisplaytest{\eq={#1}%
\if!#3!\eqnofalse\else\eqnotrue\leqnotrue\eqn={#2}\fi}
\catcode`@=12 
  \mathchardef\Gamma="0100
  \mathchardef\Delta="0101
  \mathchardef\Theta="0102
  \mathchardef\Lambda="0103
  \mathchardef\Xi="0104
  \mathchardef\Pi="0105
  \mathchardef\Sigma="0106
  \mathchardef\Upsilon="0107
  \mathchardef\Phi="0108
  \mathchardef\Psi="0109
  \mathchardef\Omega="010A

\def\sq{\hbox{\rlap{$\sqcap$}$\sqcup$}}

\def\utw{\smash{\rlap{\lower5pt\hbox{$\sim$}}}}
\def\udtw{\smash{\rlap{\lower6pt\hbox{$\approx$}}}}

\def\diameter{{\ifmmode\mathchoice
{\ooalign{\hfil\hbox{$\displaystyle/$}\hfil\crcr
{\hbox{$\displaystyle\mathchar"20D$}}}}
{\ooalign{\hfil\hbox{$\textstyle/$}\hfil\crcr
{\hbox{$\textstyle\mathchar"20D$}}}}
{\ooalign{\hfil\hbox{$\scriptstyle/$}\hfil\crcr
{\hbox{$\scriptstyle\mathchar"20D$}}}}
{\ooalign{\hfil\hbox{$\scriptscriptstyle/$}\hfil\crcr
{\hbox{$\scriptscriptstyle\mathchar"20D$}}}}
\else{\ooalign{\hfil/\hfil\crcr\mathhexbox20D}}%
\fi}}
\def\bbbr{{\rm I\!R}} 

\def\bbbn{{\rm I\!N}} 

\def\bbbone{{\mathchoice {\rm 1\mskip-4mu l} {\rm 1\mskip-4mu l}
{\rm 1\mskip-4.5mu l} {\rm 1\mskip-5mu l}}}
\def\bbbc{{\mathchoice {\setbox0=\hbox{$\displaystyle\rm C$}\hbox{\hbox
to0pt{\kern0.4\wd0\vrule height0.9\ht0\hss}\box0}}
{\setbox0=\hbox{$\textstyle\rm C$}\hbox{\hbox
to0pt{\kern0.4\wd0\vrule height0.9\ht0\hss}\box0}}
{\setbox0=\hbox{$\scriptstyle\rm C$}\hbox{\hbox
to0pt{\kern0.4\wd0\vrule height0.9\ht0\hss}\box0}}
{\setbox0=\hbox{$\scriptscriptstyle\rm C$}\hbox{\hbox
to0pt{\kern0.4\wd0\vrule height0.9\ht0\hss}\box0}}}}
\def\bbbe{{\mathchoice {\setbox0=\hbox{\smalletextfont e}\hbox{\raise
0.1\ht0\hbox to0pt{\kern0.4\wd0\vrule width0.3pt height0.7\ht0\hss}\box0}}
{\setbox0=\hbox{\smalletextfont e}\hbox{\raise
0.1\ht0\hbox to0pt{\kern0.4\wd0\vrule width0.3pt height0.7\ht0\hss}\box0}}
{\setbox0=\hbox{\smallescriptfont e}\hbox{\raise
0.1\ht0\hbox to0pt{\kern0.5\wd0\vrule width0.2pt height0.7\ht0\hss}\box0}}
{\setbox0=\hbox{\smallescriptscriptfont e}\hbox{\raise
0.1\ht0\hbox to0pt{\kern0.4\wd0\vrule width0.2pt height0.7\ht0\hss}\box0}}}}
\def\bbbq{{\mathchoice {\setbox0=\hbox{$\displaystyle\rm Q$}\hbox{\raise
0.15\ht0\hbox to0pt{\kern0.4\wd0\vrule height0.8\ht0\hss}\box0}}
{\setbox0=\hbox{$\textstyle\rm Q$}\hbox{\raise
0.15\ht0\hbox to0pt{\kern0.4\wd0\vrule height0.8\ht0\hss}\box0}}
{\setbox0=\hbox{$\scriptstyle\rm Q$}\hbox{\raise
0.15\ht0\hbox to0pt{\kern0.4\wd0\vrule height0.7\ht0\hss}\box0}}
{\setbox0=\hbox{$\scriptscriptstyle\rm Q$}\hbox{\raise
0.15\ht0\hbox to0pt{\kern0.4\wd0\vrule height0.7\ht0\hss}\box0}}}}
\def\bbbt{{\mathchoice {\setbox0=\hbox{$\displaystyle\rm
T$}\hbox{\hbox to0pt{\kern0.3\wd0\vrule height0.9\ht0\hss}\box0}}
{\setbox0=\hbox{$\textstyle\rm T$}\hbox{\hbox
to0pt{\kern0.3\wd0\vrule height0.9\ht0\hss}\box0}}
{\setbox0=\hbox{$\scriptstyle\rm T$}\hbox{\hbox
to0pt{\kern0.3\wd0\vrule height0.9\ht0\hss}\box0}}
{\setbox0=\hbox{$\scriptscriptstyle\rm T$}\hbox{\hbox
to0pt{\kern0.3\wd0\vrule height0.9\ht0\hss}\box0}}}}
\def\bbbs{{\mathchoice
{\setbox0=\hbox{$\displaystyle     \rm S$}\hbox{\raise0.5\ht0\hbox
to0pt{\kern0.35\wd0\vrule height0.45\ht0\hss}\hbox
to0pt{\kern0.55\wd0\vrule height0.5\ht0\hss}\box0}}
{\setbox0=\hbox{$\textstyle        \rm S$}\hbox{\raise0.5\ht0\hbox
to0pt{\kern0.35\wd0\vrule height0.45\ht0\hss}\hbox
to0pt{\kern0.55\wd0\vrule height0.5\ht0\hss}\box0}}
{\setbox0=\hbox{$\scriptstyle      \rm S$}\hbox{\raise0.5\ht0\hbox
to0pt{\kern0.35\wd0\vrule height0.45\ht0\hss}\raise0.05\ht0\hbox
to0pt{\kern0.5\wd0\vrule height0.45\ht0\hss}\box0}}
{\setbox0=\hbox{$\scriptscriptstyle\rm S$}\hbox{\raise0.5\ht0\hbox
to0pt{\kern0.4\wd0\vrule height0.45\ht0\hss}\raise0.05\ht0\hbox
to0pt{\kern0.55\wd0\vrule height0.45\ht0\hss}\box0}}}}
\def\bbbz{{\mathchoice {\hbox{$\sans\textstyle Z\kern-0.4em Z$}}
{\hbox{$\sans\textstyle Z\kern-0.4em Z$}}
{\hbox{$\sans\scriptstyle Z\kern-0.3em Z$}}
{\hbox{$\sans\scriptscriptstyle Z\kern-0.2em Z$}}}}
\def\qed{\ifmmode\sq\else{\unskip\nobreak\hfil
\penalty50\hskip1em\null\nobreak\hfil\sq
\parfillskip=0pt\finalhyphendemerits=0\endgraf}\fi}
\newfam\sansfam
\textfont\sansfam=\tensans\scriptfont\sansfam=\sevensans
\scriptscriptfont\sansfam=\fivesans
\def\sans{\fam\sansfam\tensans}
\def\stackfigbox{\if
Y\FIG\global\setbox\figbox=\vbox{\unvbox\figbox\box1}%
\else\global\setbox\figbox=\vbox{\box1}\global\let\FIG=Y\fi}
\def\placefigure{\global\let\twolegs=N
\dimen0=\ht1\advance\dimen0by\dp1
\advance\dimen0by5\baselineskip
\advance\dimen0by0.4true cm
\ifdim\dimen0>\vsize\pageinsert\box1\vfill\endinsert
\else
\if Y\FIG\stackfigbox\else
\dimen0=\pagetotal\ifdim\dimen0<\pagegoal
\advance\dimen0by\ht1\advance\dimen0by\dp1\advance\dimen0by1.4true cm
\ifdim\dimen0>\pagegoal\stackfigbox
\else\box1\vskip4true mm\fi
\else\box1\vskip4true mm\fi\fi\fi\egroup}
\let\twolegs=N
\def\begfig#1cm#2\endfig{\par\bgroup
\setbox1=\vbox{\kern1true cm\hbox{\vrule height#1true cm}#2}\placefigure}
\def\begdoublefig#1cm #2 #3 \enddoublefig{\begfig#1cm%
\vskip-.8333\baselineskip\line{\vtop{\hsize=6.95true cm#2}\hfill
\global\let\twolegs=N\vtop{\hsize=6.95true cm#3}}\endfig}
\def\begfigside#1cm#2cm#3\endfigside{\dimen0=#2true cm
\ifdim\dimen0<0.4\hsize\message{Room for legend to narrow;
begfigside changed to begfig}\begfig#1cm#3\endfig\else\bgroup
\par\def\figure##1##2{\vbox{\noindent\petit{\bf
Fig.\ts##1\unskip.\ }\ignorespaces ##2\par}}%
\dimen0=\hsize\advance\dimen0 by-.8true cm\advance\dimen0 by-#2true cm
\setbox1=\vbox{\kern1true cm\hbox{\hbox to\dimen0{\vrule
height#1true cm\hrulefill}%
\kern.8true cm\vbox{\hsize=#2true cm#3}}}\placefigure\fi}
\def\begfigsidetop#1cm#2cm#3\endfigsidetop{\dimen0=#2true cm
\ifdim\dimen0<0.4\hsize\message{Room for legend to narrow; begfigsidetop
changed to begfig}\begfig#1cm#3\endfig\else\bgroup
\par\def\figure##1##2{\vbox{\noindent\petit{\bf
Fig.\ts##1\unskip.\ }\ignorespaces ##2\par}}%
\dimen0=\hsize\advance\dimen0 by-.8true cm\advance\dimen0 by-#2true cm
\setbox1=\vbox{\hbox{\hbox to\dimen0{\vrule height#1true cm\hrulefill}%
\kern.8true cm\vbox to#1true cm{\hsize=#2 true cm#3\vfill
}}}\placefigure\fi}
\def\figure#1#2{\if Y\twolegs\else\vskip6pt\fi
\vbox{\noindent\petit{\bf Fig.\ts#1\unskip.\
}\ignorespaces #2\vskip\normalbaselineskip}\global\let\twolegs=Y}

\def\begtab#1cm#2\endtab{\par
\ifvoid\topins\midinsert\medskip\vbox{#2\hrule width2true cm
\kern#1true cm\hrule width2true cm}\endinsert
\else\topinsert\vbox{#2\kern#1true cm}\endinsert\fi}
\def\begpet{\vskip6pt\bgroup\petit}
\def\endpet{\vskip6pt\egroup}
\newdimen\refindent
\def\begrefchapter#1{\titleb{}{\ignorespaces#1}%
\bgroup\petit
\setbox0=\hbox{1000.\enspace}\refindent=\wd0}
\def\begrefbook#1{\titlea{}{\ignorespaces#1}%
\bgroup\petit
\setbox0=\hbox{1000.\enspace}\refindent=\wd0}

\def\ref{\goodbreak
\hangindent\oldparindent\hangafter=1
\noindent\ignorespaces}
\def\refno#1{\goodbreak
\hangindent\refindent\hangafter=1
\noindent\hbox to\refindent{#1\hss}\ignorespaces}
\def\vec#1{{\textfont1=\tams\scriptfont1=\tamss
\textfont0=\tenbfne\scriptfont0=\sevenbfne
\mathchoice{\hbox{$\displaystyle#1$}}{\hbox{$\textstyle#1$}}
{\hbox{$\scriptstyle#1$}}{\hbox{$\scriptscriptstyle#1$}}}}
\def\petit{\def\rm{\fam0\eightrm}%
\textfont0=\eightrm \scriptfont0=\sixrm \scriptscriptfont0=\fiverm
 \textfont1=\eighti \scriptfont1=\sixi \scriptscriptfont1=\fivei
 \textfont2=\eightsy \scriptfont2=\sixsy \scriptscriptfont2=\fivesy
 \def\it{\fam\itfam\eightit}%
 \textfont\itfam=\eightit
 \def\sl{\fam\slfam\eightsl}%
 \textfont\slfam=\eightsl
 \def\bf{\fam\bffam\eightbfne}%
 \textfont\bffam=\eightbfne \scriptfont\bffam=\sixbfne
 \scriptscriptfont\bffam=\fivebfne
 \def\sans{\fam\sansfam\eightsans}%
 \textfont\sansfam=\eightsans \scriptfont\sansfam=\sixsans
 \scriptscriptfont\sansfam=\fivesans
 \def\tt{\fam\ttfam\eighttt}%
 \textfont\ttfam=\eighttt
 \normalbaselineskip=10pt
 \setbox\strutbox=\hbox{\vrule height7pt depth2pt width0pt}%
 \normalbaselines\rm
\def\vec##1{{\textfont1=\tbms\scriptfont1=\tbmss
\textfont0=\eightbf\scriptfont0=\sixbf
\mathchoice{\hbox{$\displaystyle##1$}}{\hbox{$\textstyle##1$}}
{\hbox{$\scriptstyle##1$}}{\hbox{$\scriptscriptstyle##1$}}}}}
\headline={\petit\def\newline{ }\def\fonote#1{}\ifodd\pageno
\hfil\botmark\unskip\kern1.8true cm\llap{\folio}\else\leftheadline\fi}
\def\leftheadline{\rlap{\folio}\kern1.8true cm Missing TITELA\hfil}
\mark{ }
\nopagenumbers
%
\let\header=Y
\let\FIG=N
\newbox\figbox
\output={\if N\header\headline={\hfil}\fi\plainoutput\global\let\header=Y
\if Y\FIG\topinsert\unvbox\figbox\endinsert\global\let\FIG=N\fi}
\def\titlearunning#1{\message{Running head on left hand sides (TITLEA)
has been changed}\gdef\leftheadline{\rlap{\folio}\kern1.8true
cm\ignorespaces#1\hfil}{\def\newline{ }\def\fonote##1{}\mark{ #1}}%
\ignorespaces}
\def\titlebrunning#1{\message{Running head on right hand sides (TITLEB)
has been changed}\mark{#1}\ignorespaces}
\let\lasttitle=N
\def\author#1{\vfill\supereject
     \bgroup
     \baselineskip=22pt
     \lineskip=0pt
     \pretolerance=10000
     \authfont
     \rightskip 0pt plus 6em
     \noindent
     \ignorespaces#1\vskip2true cm\egroup}
\def\head#1#2{\bgroup
     \baselineskip=36pt
     \lineskip=0pt
     \pretolerance=10000
     \headfont
     \rightskip 0pt plus 6em
     \noindent
     \ignorespaces#1\vskip1true cm
     \baselineskip=22pt
     \noindent\subhfont#2\vfill
     \parindent=0pt
     \baselineskip=16pt
     \markfont Springer-Verlag\newline
     Berlin Heidelberg New York\newline
     London Paris Tokyo\par
     \egroup\let\header=N\eject}
\def\titlea#1#2{\vfill\supereject\let\header=N
     \bgroup
\textfont0=\tafontt \scriptfont0=\tafonts \scriptscriptfont0=\tafontss
\textfont1=\tamt \scriptfont1=\tams \scriptscriptfont1=\tamss
\textfont2=\tast \scriptfont2=\tass \scriptscriptfont2=\tasss
\textfont3=\tasyt \scriptfont3=\tasys \scriptscriptfont3=\tenex
     \baselineskip=20pt
     \lineskip=0pt
     \pretolerance=10000
     \tafontt
     \rightskip 0pt plus 6em
     \noindent
     \if!#1!\ignorespaces#2
     \else\setbox0=\hbox{\ignorespaces#1\unskip\
     }\hangindent=\wd0
     \hangafter=1\box0\ignorespaces#2\fi
     \vskip90pt\egroup
     \nobreak
     \parindent=0pt
     \everypar={\global\parindent=\oldparindent
     \global\let\lasttitle=N\global\everypar={}}%
     \global\let\lasttitle=A%
     \setbox0=\hbox{\petit\def\newline{ }\def\fonote##1{}\kern1.8true
     cm\ignorespaces#2}\ifdim\wd0>\hsize
     \message{Your TITLEA exceeds the headline, please use a short form with
TITLEARUNNING}\gdef\leftheadline{\rlap{\folio}\kern1.8true cm
TITLEA suppressed due to excessive length\hfil}%
\mark{TITLEA suppressed due to excessive length}%
\else{\def\newline{ }\def\fonote##1{}\mark{ #2}}\fi
     \ignorespaces}
\def\motto#1#2{\vskip-28pt\if M\lasttitle\vskip18pt\fi
\bgroup\petit\leftskip=90pt\noindent\ignorespaces#1
\if!#2!\else\medskip\noindent\it\ignorespaces#2\fi\vskip28pt\egroup
\let\lasttitle=M
\parindent=0pt
\everypar={\global\parindent=\oldparindent
\global\let\lasttitle=N\global\everypar={}}%
\global\let\lasttitle=M%
\ignorespaces}
 \def\titleb#1#2{\if N\lasttitle\else\vskip-28pt
     \fi
     \vskip23pt plus 4pt minus4pt
     \bgroup
\textfont0=\tbfontt \scriptfont0=\tbfonts \scriptscriptfont0=\tbfontss
\textfont1=\tbmt \scriptfont1=\tbms \scriptscriptfont1=\tbmss
\textfont2=\tbst \scriptfont2=\tbss \scriptscriptfont2=\tbsss
\textfont3=\tasys \scriptfont3=\tenex \scriptscriptfont3=\tenex
     \baselineskip=18pt
     \lineskip=18pt
     \pretolerance=10000
     \noindent
     \tbfontt
     \rightskip 0pt plus 6em
     \setbox0=\vbox{\vskip23pt\def\fonote##1{}%
     \noindent
     \if!#1!\ignorespaces#2
     \else\setbox0=\hbox{\ignorespaces#1\unskip\ }\hangindent=\wd0
     \hangafter=1\box0\ignorespaces#2\fi
     \vskip16pt}%
     \dimen0=\pagetotal\advance\dimen0 by-\pageshrink
     \ifdim\dimen0<\pagegoal
     \dimen0=\ht0\advance\dimen0 by\dp0\advance\dimen0 by
     3\normalbaselineskip
     \advance\dimen0 by\pagetotal
     \ifdim\dimen0>\pagegoal\vfill\eject\fi\fi
     \noindent
     \if!#1!\ignorespaces#2
     \else\setbox0=\hbox{\ignorespaces#1\unskip\ }\hangindent=\wd0
     \hangafter=1\box0\ignorespaces#2\fi
     \vskip16pt plus4pt minus4pt\egroup
     \nobreak
     \parindent=0pt
     \everypar={\global\parindent=\oldparindent
     \global\let\lasttitle=N\global\everypar={}}%
     \global\let\lasttitle=B%
     \setbox0=\hbox{\petit\def\newline{ }\def\fonote##1{}\kern1.8true
     cm\ignorespaces#2}\ifdim\wd0>\hsize
     \message{Your TITLEB exceeds the headline,  please use a short form with
TITLEBRUNNING}\mark{TITLEB suppressed due to excessive length}%
\else{\def\newline{ }\def\fonote##1{}\mark{ #2}}\fi
     \ignorespaces}
 \def\titlec#1#2{\if N\lasttitle\else\vskip-4pt\vskip-\baselineskip
     \fi
     \vskip15pt plus 4pt minus4pt
     \bgroup
     \pretolerance=10000
     \noindent
\textfont0=\tenbf \scriptfont0=\sevenbf \scriptscriptfont0=\fivebf
\textfont1=\tams \scriptfont1=\tamss \scriptscriptfont1=\tbmss
     \tenbf
     \rightskip 0pt plus 6em
     \setbox0=\vbox{\vskip 15pt\def\fonote##1{}%
     \noindent
     \if!#1!\ignorespaces#2
     \else\setbox0=\hbox{\ignorespaces#1\ }\hangindent=\wd0
     \hangafter=1\box0#2\fi
     \vskip8pt}%
     \dimen0=\pagetotal\advance\dimen0 by-\pageshrink
     \ifdim\dimen0<\pagegoal
     \dimen0=\ht0\advance\dimen0 by\dp0\advance\dimen0 by 3\baselineskip
     \advance\dimen0 by\pagetotal
     \ifdim\dimen0>\pagegoal\vfill\eject\fi\fi
     \noindent
     \if!#1!\ignorespaces#2
     \else\setbox0=\hbox{\ignorespaces#1\ }\hangindent=\wd0
     \hangafter=1\box0#2\fi
     \vskip8pt plus4pt minus4pt\egroup
     \nobreak
     \global\let\lasttitle=C%
     \parindent=0pt
     \everypar={\global\parindent=\oldparindent
     \global\let\lasttitle=N\global\everypar={}}%
     \ignorespaces}
 \def\titled#1#2{\if N\lasttitle\else\vskip-3pt\vskip-\baselineskip
     \fi
     \vskip15pt plus 4pt minus 4pt
     \bgroup
\textfont0=\tenbfne \scriptfont0=\sevenbfne \scriptscriptfont0=\fivebfne
\textfont1=\tams \scriptfont1=\tamss \scriptscriptfont1=\tbmss
     \tenbfne
     \noindent
     \if!#1!\else\ignorespaces#1\unskip\ \fi
     \ignorespaces#2\quad\egroup
     \ignorespaces}
\let\ts=\thinspace
\def\footnoterule{\kern-3pt\hrule width 2true cm\kern2.6pt}
\newcount\footcount \footcount=1
\def\advftncnt{\advance\footcount by1\global\footcount=\footcount}
\def\fonote#1{$^{\the\footcount}$\begingroup\petit
\def\line##1{{\advance\hsize by-0.5\oldparindent\hbox to\hsize{##1}}}%
\def\textindent##1{\hangindent0.5\oldparindent\noindent\hbox
to0.5\oldparindent{##1\hss}\ignorespaces}%
\vfootnote{$^{\the\footcount}$}{#1}\endgroup\advftncnt}
\def\item#1{\par\noindent
\hangindent6.5 mm\hangafter=0
\llap{#1\enspace}\ignorespaces}








\long\def\lemma#1#2{\removelastskip\vskip\baselineskip\noindent{\tenbfne
Lemma\if!#1!\else\ #1\fi.\quad}\ignorespaces#2\vskip\baselineskip}
\long\def\proposition#1#2{\removelastskip\vskip\baselineskip\noindent{\tenbfn
e
Proposition\if!#1!\else\ #1\fi.\quad}\ignorespaces#2\vskip\baselineskip}
\long\def\theorem#1#2{\removelastskip\vskip\baselineskip\noindent{\tenbfne
Theorem\if!#1!\else\ #1\fi.\quad}\ignorespaces#2\vskip\baselineskip}
\long\def\corollary#1#2{\removelastskip\vskip\baselineskip\noindent{\tenbfne
Corollary\if!#1!\else\ #1\fi.\quad}\ignorespaces#2\vskip\baselineskip}
\long\def\example#1#2{\removelastskip\vskip\baselineskip\noindent{\tenbfne
Example\if!#1!\else\ #1\fi.\quad}\ignorespaces#2\vskip\baselineskip}
\long\def\exercise#1#2{\removelastskip\vskip\baselineskip\noindent{\tenbfne
Exercise\if!#1!\else\ #1\fi.\quad}\ignorespaces#2\vskip\baselineskip}
\long\def\problem#1#2{\removelastskip\vskip\baselineskip\noindent{\tenbfne
Problem\if!#1!\else\ #1\fi.\quad}\ignorespaces#2\vskip\baselineskip}
\long\def\solution#1#2{\removelastskip\vskip\baselineskip\noindent{\tenbfne
Solution\if!#1!\else\ #1\fi.\quad}\ignorespaces#2\vskip\baselineskip}
\long\def\proof{\removelastskip\vskip\baselineskip\noindent{\it
Proof.\quad}\ignorespaces}
\long\def\remark#1#2{\removelastskip\vskip\baselineskip\noindent{\it
Remark\if!#1!\else\ #1\fi.\quad}\ignorespaces#2\vskip\baselineskip}
\long\def\definition#1#2{\removelastskip\vskip\baselineskip\noindent{\tenbfne
Definition\if!#1!\else\
#1\fi.\quad}{\it\ignorespaces#2}\vskip\baselineskip}
\def\frame#1{\bigskip\vbox{\hrule\hbox{\vrule\kern\parindent
\vbox{\kern5pt\advance\hsize by-0.8pt\advance\hsize by-2\parindent
\centerline{\vbox{\noindent#1}}\kern5pt}\kern\parindent\vrule}\hrule}\bigskip
}
\def\framedformula#1#2{$$\vcenter{\hrule\hbox{\vrule\kern5pt
\vbox{\kern5pt\hbox{$\displaystyle#1$}%
\kern5pt}\kern5pt\vrule}\hrule}\eqno#2$$}
\def\typeset{\noindent{\petit This book was processed by the author
using the \TeX\ Macropackage from Springer-Verlag.\par}}
\outer\def\byebye{\bigskip\bigskip\typeset
\footcount=1\ifx\speciali\undefined\else
\loop\smallskip\noindent special character No\number\footcount:
\csname special\romannumeral\footcount\endcsname
\advance\footcount by 1\global\footcount=\footcount
\ifnum\footcount<11\repeat\fi
\par\vfill\supereject\end}

\newcount\chapternumber \chapternumber=0
\newcount\sectionnumber \sectionnumber=0
\newcount\subsectionnumber \subsectionnumber=0
\newcount\picturenumber \picturenumber = 0
\newcount\topicnumber \topicnumber=0
\newcount\notenumber \notenumber=0
\newcount\referencenumber \referencenumber=0
\newcount\refnumber \refnumber=0
\newcount\enumerationnumber \enumerationnumber=0

\immediate\openout0=References
\immediate\openout1=Topics

\def\today{\ifcase\month \or Janvier\or Fevrier \or Mars \or Avril \or
May\or Juin \or Juli \August \or Septembre \or Octobre \or Novembre \or
Decembre\fi \space\number\today, \number\year}

\def\macro#1{\csname#1\endcsname}

\def\Numbering{\incr\topicnumber
                       \hfill{\tenbfne (\the\sectionnumber.\the\topicnumber)}}

\def\Topref#1{\incr\topicnumber
		\global\expandafter\edef\csname#1\endcsname
                {(\the\sectionnumber.\the\topicnumber)}
\immediate\write1{\def#1{(\the\sectionnumber.\the\topicnumber)}}
\decr\topicnumber}

\def\Sectionref#1{\incr\sectionnumber
		\global\expandafter\edef\csname#1\endcsname
                {\the\sectionnumber}
\immediate\write1{\def#1{\the\sectionnumber}}
     \decr\sectionnumber}

\def\Pictref#1{\incr\picturenumber
		\global\expandafter\edef\csname#1\endcsname
                {(\the\sectionnumber.\the\picturenumber)}
\immediate\write1{\def#1{(\the\sectionnumber.\the\picturenumber)}}
     \decr\picturenumber}

\def\note#1{
\incr\notenumber\footnote{$^{\the\notenumber}$}{\rm #1}}

\def\incr#1{\global\advance#1 by1}
\def\decr#1{\global\advance#1 by-1}

\def\Chapter#1{\incr\chapternumber \topicnumber = 0
\notenumber = 0 \picturenumber = 0 \sectionnumber=0 \subsectionnumber=0
\titlea{\uppercase\expandafter{\romannumeral\chapternumber}}{#1} }

\def\Section#1{\topicnumber=0 \picturenumber = 0 \subsectionnumber=0
       		\incr\sectionnumber \titleb{\the\sectionnumber.}{#1}
		\titlebrunning\papertitel }

\def\Subsection#1{\medskip{\bf #1} }

\def\Pref#1{}

\def\Titel#1{\topicnumber = 0
\notenumber = 0 \picturenumber = 0
\titlearunning{#1} }

\long\def\Exercise#1{\incr\topicnumber
\removelastskip\vskip\baselineskip\noindent{\tenbfne
(\the\sectionnumber.\the\topicnumber)
Exercise.\quad}\ignorespaces#1\vskip\baselineskip\par}

\long\def\Theorem#1{\incr\topicnumber
\medbreak\noindent{\tenbfne
(\the\sectionnumber.\the\topicnumber)
Theorem.\quad}{\it\ignorespaces  #1\vskip\baselineskip}\medbreak}

\long\def\Proposition#1{\incr\topicnumber
\medbreak\noindent{\tenbfne
(\the\sectionnumber.\the\topicnumber) Proposition.\quad}
 {\it \ignorespaces#1}\medbreak}

\long\def\Definition#1{\incr\topicnumber
\medbreak\noindent{\tenbfne
(\the\sectionnumber.\the\topicnumber)
Definition.\quad}{\it\ignorespaces#1}\medbreak}

\long\def\Lemma#1{\incr\topicnumber
\medbreak\noindent{\tenbfne
(\the\sectionnumber.\the\topicnumber)
Lemma.\quad} {\it\ignorespaces#1}\medbreak}

\long\def\Corollary#1{\incr\topicnumber
\medbreak\noindent{\tenbfne
(\the\sectionnumber.\the\topicnumber)
Corollary. \quad} {\it\ignorespaces#1}\medbreak}

\long\def\Example{\incr\topicnumber
\medbreak\noindent{\tenbfne
(\the\sectionnumber.\the\topicnumber)  Example.\quad }\ignorespaces}

\def\Proof{\proof}

\def\firstcitation#1#2{\incr\referencenumber
		 \lastbox\unskip~\hbox{[\the\referencenumber]}
                \global\expandafter\edef\csname#1\endcsname
             {\lastbox\unskip~[\the\referencenumber]}
                \immediate\write0{Ref #2 par}\hskip-1em\ }

\def\Legend#1{\incr\picturenumber
\figure{\the\sectionnumber.\the\picturenumber}{#1} }

\def\Ref#1{	\incr\refnumber \item{[\the\refnumber]} #1 \medskip}

\catcode\lq\_=8


\def\K{K} 
\def\H{{\cal W}} 

\def\Translat#1{T_{#1}} 


\def\Id{\bbbone} 


\def\U{U}
\def\R{{\bbbr}} 

\def\Z{{\bbbz}}  
\def\N{{\bbbn}}  


\def\L#1{L^{#1}}   

\def\Torus{{\bbbt}} 


\def\s{s}	   
\def\ss{r} 
\def\x{x}


\def\suchthat{:}



\def\i{i} 
\def\tt{{u}} 

\def\z{z} 



\def\inv#1{{1\over #1}}

\def\scalarproduct#1#2{\left\langle#1\mid#2\right\rangle}
\def\fou#1{\widehat{#1}} 

\def\conj#1{\overline{#1}} 

\def\abs#1{\left\vert #1 \right\vert} 


\def\norm#1{\left\Vert #1 \right\Vert}







\def\p{p}  
\def\e{e}   
\def\n{n}    
\def\m{m}    

\def\k{k}   
\def\A{A}   

\def\o{o}    
\def\O{O}  


\def\stopit{\write1{\sectionnumber=\the\sectionnumber}
             \write1{\topicnumber=\the\topicnumber}
              \write1{\notenumber=\the\notenumber}
               \write1{\referencenumber=\the\referencenumber}
                \write1{\pageno=\number\pageno}
 \end
}

\def\Shift{S}

\def\Ndim{N}

\def\H{H} 

\def\qmfcond{(1.2)}
\def\coomnm{(2.1)}
\def\givenbzagdghsjak{(2.2)}
\def\fouskdiaidufuuud{(2.3)}

\def\ssususiifuytdrs{(3.1)}
\def\qumfmatirx{(3.2)}
\def\ggvtzusiauzetq{(3.3)}
\def\xhjkslajhfas{(3.4)}
\def\satifufufhduzsiuazep{(3.6)}
\def\thehehehdjdjdjdksjhu{(3.7)}
\def\ghjftrzuzdsts{(3.8)}
\def\hfhfdjslakjhfdjkz{(5.1)}

\def\hhgtwzsastenm{(7.2)}
\def\hjkldjh{(7.3)}
\def\prcxvcjkhsuoisndjf{(7.5)}
\def\beforegisdshnans{(7.6)}

\def\Unit{{\cal U}}
\def\Shift{S}
\def\MM{M}

\def\HH{\Gamma}

\def\F{F}
\def\KK{K}
\def\Rep{R}
\def\u{u}
\def\Map{{\cal M}}
\def\Fix{\hbox{Fix}}

\def\SU{\hbox{SU}}

\def\x{x}

\def\U{U}

\def\Unit{{\cal U}}

\def\A{A}
\def\B{B}
\def\CC{C}

\def\Shift{S}

\def\kk{{k^\prime}}
\def\Map{\hbox{Map}}

\let\header=N

{\headfont
\centerline{Quadrature Mirror Filters}
\bigskip
\centerline{and}
\bigskip
\centerline{Loop Groups}}
\bigskip
\bigskip {\authfont
\centerline{Matthias Holschneider\note{On leave from CNRS Centre
de Physique Th\'eorique, Luminy, Case 907, F-13288 Marseille,
e-mail:
 hols.at.marcptnx1.univ-mrs.fr}}
\bigskip
\centerline{and}
\bigskip
\centerline{Ulrich Pinkall\note{Technische Universit\"at Berlin,
FB Mathematik, Strasse des 17. Juni 135, 10623 Berlin}} }
\bigskip
\centerline{Max-Plank-Arbeitgruppe Mathematik}
\centerline{Universit\"at Potsdam}
\centerline{am Neuen Palais}
\centerline{Potsdam}

\vskip 2 cm

\def\HH{\Gamma}

{\bf Abstract:} In this paper we want to show, that the finite
impulse response  quadratic mirror filters (QMF) associated to a
tower of grids $\Gamma\subset\H=\Z^\n$  can be identified with a
right  coset of the subgroup
$\Fix(\Translat{\HH^\perp},\Map(\Torus^\n\to\U(\Ndim)\suchthat\hbox{poly})$
of the group of polynomial loops
$\Map(\Torus^\n\to\U(\Ndim)\suchthat\hbox{poly})$  with
$\Ndim=\abs{\H/\Gamma}$.  The QMF with some vanishing moments can
be identified with cosets of subgroups.  The problem to
parameterize all finite impulse response QMF in arbitrary space
dimensions is now equivalent to factorize all polynomial loops.
\vfill
\eject

\let\header=Y
\pageno=1

\Titel{Quadratic Mirror Filters and Loop Groups}
\def\papertitel{Matthias Holschneider and Ulrich Pinkall}

\Section{Introduction: the case $\Z/2\Z\subset\Z$}

The kind of theorem  we propose to proof  in this paper is  the
easiest illustrated with the one dimensional case.
 We inclose it here for the sake of clarity although all results
are  well known. After a short introduction to QMF filters we want
to show that the  construction of all QMF is equivalent to
construct the  group of unitary operators that commute with the
translations. This in turn will be equivalent to the construction
of loops in $\U(2)$.

Consider the space of square integrable sequences $\L2(\Z)$. At
the heart of QMF-filters  \def\Gastand{[1]}\Gastand\  and discrete
wavelet analysis (e.g.  \def\Ingrid{[2]}\Ingrid) is the idea to
split  any sequence into an even and odd part; that is we write
$$
\L2(\Z) = \L2(2\Z) \oplus \L2(2\Z+1).
$$ The associated orthogonal projectors on $\L2(2\Z)$ and on
$\L2(2\Z+1)$ shall be denoted by $\Pi_0$ and
$\Pi_1=\Translat1\Pi_0\Translat1^\ast$ where
$\Translat\n$ is the translation operator.  Clearly no information
on a sequence $\s$ is lost if we  split it into its even and odd
sub-sequences. In particular the energy
$\norm{\s}^2$ is just the sum of the even and the odd part
$$
\norm{\s}^2 = \norm{\Pi_0\s}^2 + \norm{\Pi_1\s}^2.
$$ We now ask the question whether we can find a pre-treatment of
$\s$, such that again no information is lost and such that again
the energy is conserved.
 To be more specific we are looking for sequences $\eta_0$ and
$\eta_1$ such that the split of the filtered sequences
$\eta_0\ast\s$ and
$\eta_1\ast\s$ is again energy conserving. Thus we want to have
$$
\norm{\s}^2 = \norm{\Pi_0(\eta_0\ast\s)}^2 +
\norm{\Pi_1(\eta_1\ast\s)}^2,
\Topref{spliththjk}\Numbering
$$ with the convolution product defined as
$$
\ss\ \ast\ \s(\n) = \sum_{\k\in\Z}\ss(\n-\k) \,\s(\k) =
\scalarproduct{\Translat\n\tilde\ss}{\s}\quad\hbox{with}\ \
\tilde\ss(\k) = \conj\ss(-\k)
$$ But this requirement precisely means that we should have
$$\eqalign{
\norm{\s}^2 &= \sum_{\k\in\Z}\abs{\s(\k)}^2
\cr &= \sum_{\k\in2\Z} \abs{\eta_0\ast\s(\k)}^2 +
\sum_{\k\in2\Z+1} \abs{\eta_1\ast\s(\k)}^2
\cr &= \sum_{\k\in2\Z}
\abs{\scalarproduct{\Translat\k\tilde\eta_0}{\s}}^2 +
\sum_{\k\in2\Z+1}
\abs{\scalarproduct{\Translat\k\tilde\eta_1}{\s}}^2. }
$$ Therefore if we set
$$
\phi_0 = \tilde\eta_0, \ \phi_1 = \Translat1^\ast\tilde\eta_1,
$$ then this precisely means that
 the set of $2\Z$-translated functions satisfies at
$$
\{\Translat\k\phi_0\suchthat\k\in2\Z\}\cup\{\Translat\k\phi_1
\suchthat \k\in2\Z\}\quad\hbox{is an o.n.b. of
$\L2(\Z)$}\Topref{qmfcond}\Numbering
$$

This gives rise to a definition

\Definition{A pair of sequences in $\L2(\Z)$ satisfying \qmfcond\
is called a QMF-system (=Quadratic Mirror Filter system).}

The reason for the name \lq\lq Quadratic Mirror Filter\rq\rq\
becomes clear
 later on when we rewrite the QMF condition in Fourier space.

\def\ii{{i^\prime}} By definition a pair of sequences $\phi_0$,
$\phi_1\in\L2(\Z)$ is a QMF system if and only if
$$
\eqalign{i&)\ \hbox{the set}\  \{\Translat{2\k}\phi_0,
\Translat{2\k}\phi_1
\suchthat \k\in\Z\}\ \hbox{is complet}\cr ii&)\
\scalarproduct{\Translat\k\phi_\i}{\Translat\kk\phi_\ii} =
\delta(\i-\ii)\,
\delta(\k-\kk) \quad\hbox{for}\quad
\i,\ii\in\{0,1\} \ \hbox{and}\ \k,\kk\in\Z.}
$$ However as we shall see later in theorem \thehehehdjdjdjdksjhu\
condition $ii$ is sufficient.

\Section{The group of unitary operators with $[\U,\Translat2]=0$.}
We now want to show that the set of QMF-systems has a natural
group structure. It turns out to be isomorph to a subgroup  of
unitary operators acting in the Hilbert space
$\L2(\Z)$. We shall denote the group of all unitary operators
acting in $\L2(\Z)$ by
$\Unit(\L2(\Z))$ Let $\U$ be a unitary operator acting on
sequences in $\L2(\Z)$. Suppose further that it  commutes with the
translations by $2$
$$ [\Translat2, \U] = \Translat2 \U - \U \Translat2 =0,
\quad\hbox{or equivalently}\quad\Translat2\,\U\,\Translat2^\ast=\U
\Topref{coomnm}\Numbering
$$ Because of the second equation this class of operators shall be
denoted by $\Fix(\Translat2,\Unit(\L2(\Z)))$. Let now $\phi_0$,
$\phi_1\in\L2(\Z)$  be  a QMF-system.  Then the image of $\phi_0$
and $\phi_1$ under any unitary $\U$ satisfying
\coomnm\ clearly is again a QMF-system. Indeed as unitary images
they are again an ortho-normal base and in addition  we  have for
$\i,\ii\in\{0,1\}$ and $\k,\kk\in\Z$
$$\eqalign{
\scalarproduct{\Translat{2\k}\U\phi_\i}{\Translat{2\kk}\U\phi_{\i^\prime}}
&=
\scalarproduct{\U\Translat{2\k}\phi_\i}{\U\Translat{2\kk}\phi_{\i^\prime}}
\cr &=
 \scalarproduct{\Translat{2\k}\phi_\i}{\Translat{2\kk}\phi_{\ii}}
 = \delta(\i-{\ii})\delta(\k-\kk) }
$$ Now clearly the sequences $\phi_0=\delta$, $\phi_1 =
\Translat1\delta$  are a QMF-system.  We claim that any other
QMF-system can be obtained as image of $\delta$ and
$\Translat1\delta$ under some unitary operator $\U$ that commutes
with the $2\Z$-translates.  To see this let $\phi_0$, $\phi_1$ be a
QMF-system. Clearly any operator that commutes with the
translations with respect to $2\Z$ is completely determined by its
image on one single fundamental domain $\Z/2\Z=\{0,1\}$.
Therefore if we set
$$
\U\delta = \phi_0,\quad \U\Translat1\delta =
\phi_1,\Topref{givenbzagdghsjak}
\Numbering
$$ then there is at most  one such unitary $\U$ acting on all
 of $\L2(\Z)$ and satisfying $[\Translat2,\U]=0$. There actually
is one and it is given by
$$
\U : \s\mapsto (\Pi_0\s)\ast\phi_0 +
(\Pi_1\s)\ast(\Translat1^\ast\phi_1).
\Topref{fouskdiaidufuuud}\Numbering
$$ It obviously commutes with $\Translat2$ and we only
 have to check its unitarity. But since it sends the orthonormal
base $\{\Translat{2\k}\delta,\Translat{2\k+1}\delta\}$ into  the
orthonormal base
$\{ \Translat{2\k}\phi_0, \Translat{2\k+1}\phi_1\}$, we are done.
Therefore we have shown:

\Topref{jjfudiosidufu}
\Theorem{There is a bijection between the group
$\Fix(\Translat2, \Unit(\L2(\Z)))$  and the QMF-systems $\phi_0$,
$\phi_1$.  This relation  is explicitly given by
\givenbzagdghsjak\ and \fouskdiaidufuuud.}

This shows that the QMF have a natural group structure. In
particular one can  produce new ones from old ones by
multiplication of their associated unitary operators.

\Section{The Fourier  space picture} The definition of Fourier
transform that we use is the following
$$
\fou\cdot\ :\  \L2(\Z) \mapsto \L2(\Torus),\quad
\fou\s(\omega) = \sum_{\k\in\Z}\s(\k) \e^{-\i\omega\k}.
$$ Now the action of $\Pi_0$ and $\Pi_1$ in Fourier space is easily
computed yielding
$$
\Pi_0 : \fou\s(\omega) \mapsto {1\over2}(\fou\s(\omega) +
\fou\s(\omega+\pi)),\quad
\Pi_1 : \fou\s(\omega) \mapsto {1\over2}(\fou\s(\omega) -
\fou\s(\omega+\pi)).
\Topref{ssususiifuytdrs}\Numbering$$ Now ortho-normality of the
$2\Z$-translates of a sequence $\phi_0$ is equivalent to
$$ (\tilde\phi_0\ \ast\ \phi_0)(2\k) = \delta(\k),
\quad\hbox{with}\quad\tilde\phi_0(\n) =\conj\phi_0(-\n)
$$ which also reads
$$
\Pi_0\ (\tilde\phi_0\ \ast\ \phi_0) = \delta.
$$ We therefore see from \ssususiifuytdrs\ that being a QMF-filter
system
\qmfcond\ implies that
$$
\abs{\fou\phi_0(\omega)}^2 + \abs{\fou\phi_0(\omega+\pi)}^2 =
2,\quad
\hbox{for almost all $\omega$}.
$$ This explains the name Quadratic Mirror Filter. In the same way
we obtain
$$
\abs{\fou\phi_1(\omega)}^2 + \abs{\fou\phi_1(\omega+\pi)}^2 = 2,
\quad
\hbox{for almost all $\omega$}.
$$ Now mutual orthogonality of the respective translates
$\Translat{2\k}\phi_0$ and
$\Translat{2\kk}\phi_1$ is equivalent to
$$
\Pi_0\ (\tilde\phi_0\ \ast\ \phi_1) = 0
$$ or in Fourier space
$$
\conj{\fou\phi}_0(\omega)\fou\phi_1(\omega) +
\conj{\fou\phi}_0(\omega+\pi)\fou\phi_1(\omega+\pi) = 0,
\quad
\hbox{for almost all $\omega$.}
$$ All these relations can be summarized by saying that the matrix
$$
\inv{\sqrt2}\
\left[\matrix{\fou\phi_0(\omega)&\fou\phi_0(\omega+\pi)\cr
	\fou\phi_1(\omega)&
\fou\phi_1(\omega+\pi)\cr}\right]\Topref{qumfmatirx}
\Numbering
$$ is unitary for all $\omega$.   As we shall see in a while
(theorem \thehehehdjdjdjdksjhu)  the inverse is also true: any
such matrix valued function $\chi:\Torus\to\U(2)$ satisfying
$$
\chi(\omega+\pi) = \chi(\omega) \left[\matrix{0&1\cr1&0\cr}\right]
\Topref{ggvtzusiauzetq}\Numbering
$$ gives rise to a QMF-system. We shall call such a matrix a
QMF-matrix.  By what we have seen, such a matrix
 defines a pair of  sequences such that their $2\Z$ translates are
an ortho-normal set. However the completeness still has to be
checked.

Since we have already established the equivalence between
QMF-systems and unitary operators $\U$ that commute with
$\Translat2$ let us now exploit this relation by writing the
action of the unitary operator $\U$ in
 Fourier space. From
\fouskdiaidufuuud\ it follows that  we have
$$
\U : \fou\s(\omega) \mapsto  {1\over2}(\fou\s(\omega) +
\fou\s(\omega+\pi))\ \fou\phi_0(\omega)  +
{1\over2}(\fou\s(\omega) - \fou\s(\omega+\pi))\
\fou\phi_1(\omega)\,
\e^{\i\omega}
$$ But now the action of $\U$ on $\fou\s(\omega+\pi)$ looks
similar and therefore if we introduce the vector
$\left[{\fou\s(\omega)\atop\fou\s(\omega+\pi)}\right]$ we
 may write
$$\eqalign{
\U\ :\ &\left[{\fou\s(\omega)\atop\fou\s(\omega+\pi)}\right]
\mapsto \cr  &{1\over2}\left[\matrix{\fou\phi_0(\omega)+
\fou\phi_1(\omega)\e^{\i\omega}&
\fou\phi_0(\omega)-\fou\phi_1(\omega)\e^{\i\omega}\cr
\fou\phi_0(\omega+\pi)-\fou\phi_1(\omega+\pi)\e^{\i\omega}&
\fou\phi_0(\omega+\pi)-\fou\phi_1(\omega+\pi)\e^{\i\omega}\cr}\right]\
\left[{\fou\s(\omega)\atop\fou\s(\omega+\pi)}\right]}
$$ or equivalently
$$
\U\ :\ \left[{\fou\s(\omega)\atop\fou\s(\omega+\pi)}\right]
\mapsto
{1\over2}\left[\matrix{\fou\phi_0(\omega)&\fou\phi_1(\omega)\cr
	\fou\phi_0(\omega+\pi)& \fou\phi_1(\omega+\pi)\cr}\right]\
\left[\matrix{1& 1\cr
	\e^{\i\omega}&-\e^{\i\omega}\cr}\right] \ \left[{\fou\s(\omega)
\atop\fou\s(\omega+\pi)}\right]
$$ This is a product of two unitary matrices and hence unitary
itself.It is a
 matrix valued function $\eta:\Torus\to\U(2)$ that satisfies at
$$
\eta(\omega+\pi) =  \left[\matrix{0&1\cr1&0\cr}\right]\
\eta(\omega)\
 \left[\matrix{0&1\cr1&0\cr}\right].
$$ The space of such mappings forms a group under pointwise
multiplication.  Vice versa suppose we are given such a
matrix-valued function.  Upon setting
$$
\U : \ \left[{\fou\s(\omega)\atop\fou\s(\omega+\pi)}\right]
\mapsto \eta(\omega)\
\left[{\fou\s(\omega)\atop\fou\s(\omega+\pi)}\right]
\Topref{xhjkslajhfas}\Numbering
$$ this defines an operator acting on $\fou\s(\omega)$ and the
above condition  precisely assures that this operator is well
defined as can be seen upon replacing
$\omega$ with $\omega+\pi$. It clearly is unitary. And the
translation operator
$\Translat2$ becomes diagonal
$$
\Translat2 : \left[{\fou\s(\omega)\atop\fou\s(\omega+\pi)}\right]
\mapsto
\left[\matrix{\e^{-2\i\pi\omega}&0\cr0&\e^{-2\i\pi\omega}\cr}\right]\,
\left[{\fou\s(\omega)\atop\fou\s(\omega+\pi)}\right]
$$ and therefore $\U$ commutes with $\Translat2$.  We therefore
have shown

\Theorem{There is a bijection of the subgroup of all unitary
operators acting in
$\L2(\Z)$ with $[\Translat2,\U]=0$ and the group of measurable
mappings $\eta : \Torus\to\U(2)$ satisfying (almost everywhere)
$$
\eta(\omega+\pi) = \left[\matrix{0&1\cr1&0\cr}\right]\
\eta(\omega)\
\left[\matrix{0&1\cr1&0\cr}\right].\Topref{satifufufhduzsiuazep}\Numbering
$$ The correspondence is given by \xhjkslajhfas.}

Suppose now that $\chi : \Torus\to\U(2)$  satisfies
\ggvtzusiauzetq. Since
$$
\rho(\omega) = {1\over\sqrt2}\left[\matrix{1& \e^{\i\omega}\cr
	1&-\e^{\i\omega}\cr}\right]
$$
 satisfies at
$$
\rho(\omega+\pi) = \left[\matrix{0&1\cr1&0\cr}\right]\
\rho(\omega),
$$ it follows that the combination $\xi(\omega) =
\rho(\omega)\,\chi(\omega)$ satisfies \satifufufhduzsiuazep. It
therefore gives rise to a unitary operator which implies that
$\chi(\omega)$ is coming from a QMF-system. To summarize  we have
shown

\Topref{thehehehdjdjdjdksjhu}\Theorem{There is a bijection
 between QMF-systems and measurable functions $\chi:
\Torus\to\U(2)$ satisfying (almost everywhere)
$$
\chi(\omega+\pi) = \chi(\omega) \left[\matrix{0&1\cr1&0\cr}\right]
\Topref{ghjftrzuzdsts}\Numbering
$$ it reads explicitly
$$
\left[\matrix{\fou\phi_0(\omega)&\fou\phi_0(\omega+\pi)\cr
	\fou\phi_1(\omega)& \fou\phi_1(\omega+\pi)\cr}\right]
 =
\chi(\omega) \left[\matrix{1&0\cr0&0\cr}\right]  +
\chi(\omega+\pi) \left[\matrix{0&0\cr0&1\cr}\right]
$$ }

\Section{QMF and loop groups.} A  map from the torus $\Torus$ to
the  space of unitary matrices $\U(2)$ is called a closed loop in
$\U(2)$.  The space of all  loops, $\Map(\Torus\to\U(2))$,  is a
group under pointwise multiplication (e.g.
\def\pppdiodufdi{[5]}\pppdiodufdi).

If we identify $\Torus$ with $\R/2\pi\Z$ then the loop condition
reads
$$
\xi(\omega)\in\U(2),\quad
\xi(\omega+2\pi) = \xi(\omega)
$$ The set of functions $\xi:\Torus\to\U(2)$ such that
$\Translat\pi\xi(\omega)=\xi(\omega+\pi)=\xi(\omega)$ form a (non
normal) subgroup that we shall denote by
$\Fix(\Translat\pi, \Map(\Torus\to\U(2)))$. This group may be
identified with
$\Map(\Torus\to\U(2))$ itself upon replacing $\xi(\omega)$ with
$\xi(\omega/2)$. Next consider the matrices of the special form
\qumfmatirx, or what is the same a loop $\chi(\omega)$ satisfying
\ggvtzusiauzetq. The set of such loops does not form a group.
However this set of loops is  a coset; to be more precise let
$\chi_1(\omega)$ and $\chi_2(\omega)$ be two  such loops. Then
$\xi(\omega) = \chi_1(\omega)\, \chi_2^\ast(\omega)$ satisfies at
$$
\xi(\omega+\pi) = \chi_1(\omega+\pi)\, \chi_2^\ast(\omega+\pi) =
\chi_1(\omega)\, \chi_2^\ast(\omega) = \xi(\omega),
$$ and is therefore an element of $\Fix(\Translat\pi,
\Map(\Torus\to\U(2)))$. In addition all loops in this subgroup can
be obtained that way since for
$\xi\in\Fix(\Translat\pi, \Map(\Torus\to\U(2)))$ we may set
$$
\chi_2(\omega) = \xi^\ast(\omega)\,\chi_1(\omega),
$$ and this loop satisfies at \ghjftrzuzdsts\ and it gives back
$\xi$ by setting
$\xi(\omega)=\chi_1(\omega)\, \chi_2^\ast(\omega)$.   To summarize
we have shown

\Theorem{The set of QMF-loops \qumfmatirx\ is a right coset of
$\Fix(\Translat\pi, \Map(\Torus\to\U(2)))$. Therefore let
$\chi(\omega)$ be any QMF-matrix \qumfmatirx. Then all others can
be obtained in a unique way as
$$
\eta(\omega) = \chi(\omega)\ \rho(\omega)\quad\hbox{where $\rho$
runns through}
\ \Fix(\Translat\pi, \Map(\Torus\to\U(2))).
$$}

Let us now look at the other combination; that is consider the set
of loops $\xi(\omega)=\chi_1^\ast(\omega)\,\chi_2(\omega)$. Note
that this set is not just the conjugate set of the previous one
since $\chi^\ast(\omega)$ does not satisfy \ggvtzusiauzetq\ any
more. We rather have now
$$
\xi(\omega+\pi) =  \left[\matrix{0&1\cr1&0\cr}\right]\
\xi(\omega)\
 \left[\matrix{0&1\cr1&0\cr}\right].
$$ This is again a (non normal) subgroup of all loops.  Its
elements are called   twisted loops. If we introduce the operator
 $\sigma$ acting on the loop group via conjugation
$$
\sigma : \rho(\omega) \mapsto \left[\matrix{0&1\cr1&0\cr}\right]
\ \rho(\omega)\ \left[\matrix{0&1\cr1&0\cr}\right]
$$ then we can write for the group of twisted loops
$\Fix(\sigma\,\Translat\pi,
\Map(\Torus\to\U(2)))$. Again we have that the set of QMF is a
coset which this time is a left-coset.

\Theorem{The set of QMF-loops \qumfmatirx\ is a left coset of
$\Fix(\sigma\,\Translat\pi, \Map(\Torus\to\U(2)))$. Therefore let
$\chi(\omega)$ be any QMF-matrix \qumfmatirx. Then all other
QMF-matrices
 can be obtained in a unique way as
$$
\eta(\omega) = \rho(\omega)\ \chi(\omega)\quad\hbox{where $\rho$
runns through}
\ \Fix(\sigma\,\Translat\pi, \Map(\Torus\to\U(2))).
$$}

\Section{Some subclasses of QMF.}  We now shall consider various
subgroups of the whole loop-group and their associated
QMF-systems.  The correspondence is made in the spirit of the
preceding theorems; that is the set of QMF matrices will  be
identified with a right coset of a given  subgroup of all loops.

{\bf a) The smooth loops.} Here smoothness means  that every
matrix element of $\xi(\omega)$ is a function in
$\CC^\infty(\Torus)$.
 The space of all smooth loops,
$\Map(\Torus\to\U(2)\suchthat\CC^\infty)$,  is a Lie group under
pointwise multiplication(e.g.
\pppdiodufdi).  As we know already these loops correspond to
QMF-systems $\phi_0$, $\phi_1$,  whose Fourier transforms are in
$\CC^\infty(\Torus)$, or what is the same, whose sequences are
arbitrary well localized in the sense that
$$
\sup_{\n\in\Z}(1+\abs\n)^\alpha\
\abs{\phi_\i(\n)}<\infty,\quad\hbox{for all
$\alpha>0$},
\quad \i\in\{0,1\}.
$$ Since
$$
\omega\to{1\over\sqrt2}\left[\matrix{1& \e^{\i\omega}\cr
	1&-\e^{\i\omega}\cr}\right]\
$$ is a smooth loop that is at the same time a QMF loop, we see
that all highly localized QMF-systems are obtained as a right
coset
$$
\omega\to{1\over\sqrt2}\left[\matrix{1& \e^{\i\omega}\cr
	1&-\e^{\i\omega}\cr}\right]\,\eta(\omega)
$$ where $\eta$ runs through
 $\Fix(\Translat\pi,\Map(\Torus\to\U(2))\suchthat\CC^\infty)$.
This last group is isomorphic to the group of smooth loops
 $\Map(\Torus\to\U(2))\suchthat\CC^\infty)$.

{\bf b.) The polynomial loops.}  This is the class of  loops for
which the matrix elements of $\xi(\omega)$ are trigonometric
 polynomials (= Laurent polynomials in $\z=\e^{\i\omega}$).  We
shall denote this set by
$\Map(\Torus\to\U(2)\suchthat\hbox{poly})$. The set of
QMF-systems with finite  impulse response (= their convolution
with a delta sequence has only a finite  number of non vanishing
terms) is a right-coset of
$\Fix(\Translat\pi,\Map(\Torus\to\U(2)\suchthat\hbox{poly}))$.
This last group is isomorphic to the   group of all polynomial
loops
$\Map(\Torus\to\U(2)\suchthat\hbox{poly})$.

Since
$$
\omega\to{1\over\sqrt2}\left[\matrix{1& \e^{\i\omega}\cr
	1&-\e^{\i\omega}\cr}\right]
$$ is a polynomial loop that defines a QMF, only the last
statement  has to be shown yet.
 To see this note that
 the subgroup of polynomial loops in $\Fix(\Translat\pi,
\Map(\Torus\to\U(2))\suchthat \hbox{poly})$ is isomorphically
mapped onto the group of all polynomial loops via the
identification
$\eta(\omega)\to\eta(\omega/2)$. Indeed from
$\eta(\omega+\pi)=\eta(\omega)$ is follows easily that
$\eta(\omega)$ must be a polynomial in $\e^{2\i\pi\omega}$.

{\bf c.) The vanishing moment loops.} Consider the subset of
smooth loops that satisfy
$$
\xi(\omega) = \Id + \o(\omega^\n)\quad(\omega\to0).
\Topref{hfhfdjslakjhfdjkz}\Numbering
$$ This is again a normal subgroup of
$\Map(\Torus\to\U(2)\suchthat\CC^\infty)$.  On the other hand
consider the class of highly localized QMF-systems such that
$$
\inv{\sqrt2}\,
\left[\matrix{\fou\phi_0(\omega)&\fou\phi_0(\omega+\pi)\cr
	\fou\phi_1(\omega)& \fou\phi_1(\omega+\pi)\cr}\right] =
\Id+\o(\omega^\n)\quad(\omega\to0),\Topref{fhfhhfzdstzdus}\Numbering
$$ or what is the same that
$$
\fou\phi_0(\omega) = \sqrt2 + \o(\omega^\n),\quad
\fou\phi_1(\omega) = \o(\omega^\n)\quad (\omega\to0).
$$ The second condition is equivalent to the fact that the first
$\n$ moments of $\phi_1$ vanish
$$
\sum_{\m\in\N}\, \m^\p\, \phi_1(\m) = 0,\quad\p=0,1,\dots\n.
$$ Again the subclass of such QMF-matrices is a right coset of the
subgroup of  smooth loops in
$\Fix(\Translat\pi,\Map(\Torus\to\U(2))\suchthat\CC^\infty)$
satisfying  \hfhfdjslakjhfdjkz.
 The analogue statement holds for the polynomial loops and the
finite
 impulse response QMF systems.

\Section{The factorization problem.}

Because of its particular importance in signal processing we shall
have a closer look at the polynomial loops. In particular we shall
give a technique how to generate all finite impulse response
QMF-systems.

Given that the polynomial
 loops are an infinite dimensional Lie group, the question arises
whether or not there is a small family of loops that generate all
polynomial loops. Equivalently the problem is how to generate in
an economic way all finite impulse response QMF. Such a
factorization exists as can be seen by general theorems about loop
groups (e.g. in \pppdiodufdi).  Instead of reproducing the
demonstration given there we present a different,  more geometric
construction.

As we have seen we may equally well consider the group of unitary
operators
 that commute with the translates with respect to the sub-lattice
$2\Z$.

Consider the  fundamental domain $\HH=\{0,1\}$ of $2\Z$ in $\Z$.
Recall that a  fundamental domain  is a set $\HH\subset\Z$ whose
translates with respect to $2\Z$ are a mutually disjoint cover of
$\Z$. Suppose that $\U$  is a unitary operator acting in $\L2(\Z)$
that commutes with $\Translat2$.  Clearly by translation
invariance such an operator is completely determined by its action
on the fundamental domain $\HH$. We now pick a unitary 2 by 2
matrix
$\u$ in $\U(2)$. It acts in a natural way in the two-dimensional
vector space
$\L2(\HH)$ which we consider in a natural way as a subspace of
$\L2(\Z)$.
 Now consider the well defined operator $\U_\u$, with
$[\U_\u,\Translat2]=0$,   whose restriction to the fundamental
domain $\HH$ is given by $\u$. Since this operator preserves the
$\L2$-norm in each of the translated  fundamental domains we find
that $\U_\u$ is unitary. We therefore have found a family of
 unitary operators that commute with the translations by $2$. This
family is indexed by
 points in the $4$-dimensional manifold $\U(2)$.

This family together with the translation operators $\Translat1$
actually generates all finite impulse response QMF.

\Theorem{Let $\U:\L2(\Z)\to\L2(\Z)$ be a  unitary operator with
$[\U, \Translat2]=0$. Suppose in addition that the associated  QMF
$\phi_0 = \U\,\delta_0$ and $\phi_1 = \U\,\delta_1$ is of finite
impulse response with length not larger than $2\Ndim$ (=each
sequence has at most
 $2\Ndim$
 non-zero coefficients). Then $\U$ can be factorized as follows
$$
\U = \Translat\p\,\U_{\u_{2\Ndim}}\,
\Translat1\,\U_{\u_{2\Ndim-1}} \dots
\, \Translat1\,\U_{\u_1} \quad \hbox{or}
$$ with some $\p\in\Z$. If all coefficients are real valued  then
the $\u_\n$ can be chosen in
$\O(2)$.}

\Proof Let $\phi_0 = \U\,\delta_0$ and $\phi_1 = \U\,\delta_1$ be
the associated
 QMF system. By an overall translation we may suppose that the
support  of $\phi_0$ and $\phi_1$ is contained in $[0,2\Ndim-1]$.
By the QMF-property  we have
$$\eqalign{
\scalarproduct{\phi_0}{\Translat{2(\Ndim-1)}\phi_0}  &=
\scalarproduct{\phi_1}{\Translat{2(\Ndim-1)}\phi_1} =
\scalarproduct{\phi_0}{\Translat{2(\Ndim-1)}\phi_1}
\cr &= \scalarproduct{\phi_1}{\Translat{2(\Ndim-1)}\phi_0} = 0}
$$ unless $\Ndim=1$. This shows that the four vectors
$$\eqalign{ &\e_0^l=[\phi_0(0),\phi_0(1)],\
\e_0^r=[\phi_0(2\Ndim-2),\phi_0(2\Ndim-1)]\cr
 &\e_1^l=[\phi_1(0),\phi_1(1)],
\ \e_1^r=[\phi_1(2\Ndim-2),\phi_1(2\Ndim-1)] }
$$  are pairwise orthogonal. We thus can find a two by two matrix
$\u_\Ndim\in\SU(2)$ such that
$$
\scalarproduct{\u_\Ndim\e_0^l}{[1,0]} =
\scalarproduct{\u_\Ndim\e_1^l}{[1,0]}
=\scalarproduct{\u_\Ndim\e_0^r}{[0,1]} =
\scalarproduct{\u_\Ndim\e_1^l}{[0,1]} = 0.
$$ Therefore $\U_{\u_\Ndim}\phi_0$ and $\U_{\u_\Ndim}\phi_1$ are
supported by a subset of
$\{1,2\Ndim-2\}$. We therefore may repeat the procedure with
$\Translat1^\ast\U_{\u_\Ndim}\phi_0$ and
$\Translat1^\ast\U_{\u_\Ndim}\phi_1$ and so on until $\Ndim=1$.
Here we are left  with a two dimensional vector space with two
orthogonal vectors.  We may again find a $\u_1\in\U(2)$ that maps
these two vectors to $[1,0]$ and $[0,1]$. Upon reverting all
operations we have found the desired factorization. The case where
all sequences are real valued follows the same way.
\qed

We note that in the one dimensional case a different factorization
has  been obtained before  in
 \def\Pollen{[4]}\Pollen.
\def\HH{\Gamma}
\Section{The general case}

We now want to formulate the results of the preceding  section for
general lattices;  by a lattice we mean  a subgroup of $\Z^\n$
that is isomorph to $\Z^\n$.

Let $\HH$ be a lattice in  $\H=\Z^\n$.  Let
$\Ndim=\abs{\H/\Gamma}$ be  the number of  elements in the
quotient group.  As before we define

\Definition{A system of
$\Ndim$ sequences $\phi_\k\in\L2(\H)$, $\k=0,1,\dots,\Ndim-1$ is
called a
 QMF-system if and only if the set of all $\HH$-translates of
these sequences is an ortho-normal basis of the whole space
$\L2(\H)$.  }

More explicitly we want to have
$$
\{\Translat\m\phi_\k \suchthat \k=0,1,\dots,\Ndim-1,
\m\in\Gamma\}\quad
\hbox{is an ortho-normal basis of }\L2(\H).
$$ which is again equivalent to
$$
\eqalign{i&)\
\Translat\m\phi_\k \suchthat \k=0,1,\dots,\Ndim-1, \m\in\Gamma
\quad\hbox{is complete}\cr ii&)\
\scalarproduct{\Translat\m\phi_\k}{\Translat{\m^\prime}\phi_{\k^\prime}}
=
\delta(\m-\m^\prime)\,\delta(\k-\k^\prime).}
$$ As we shall see again condition $ii$ alone is sufficient.

Consider now the group of unitary operators  acting  in the
Hilbert space $\L2(\H)$ that commute with the   translates with
respect to the sub-grid $\HH\subset\H$
$$ [\K,\Translat\m] = \K\, \Translat\m -
\Translat\m\,\K\quad\hbox{for all }\k\in\Gamma
$$ with the translation operator defined as
$$
\Translat\m : \L2(\H) \to \L2(\H),\quad\s(\cdot)
\mapsto\s(\cdot-\m).
$$ Let $\{\phi_\k\}$ be a QMF-system. We claim that the image of a
QMF-system  under a unitary operator $\K$ that commutes with the
$\Gamma$-translates is again a QMF-system. Indeed as image of an
orthonormal basis it is again an ortho-normal basis and in
addition we have as before
$$
\eqalign{
\scalarproduct{\Translat\k\K\phi_\i}{\Translat\kk\K\phi_\ii}&=
\scalarproduct{\K\Translat\k\phi_\i}{\K\Translat\kk\phi_\ii}\cr &=
\scalarproduct{\Translat\k\phi_\i}{\Translat\kk\phi_\ii}
=\delta(\i-\ii)\,
\delta(\k-\kk)}
$$ Now consider the quotient $\H/\HH$ that we may identify with a
fundamental domain for the action of $\Gamma$ on $\H$. It contains
$\Ndim=\abs{\H/\Gamma}$ points. Clearly the family of delta
sequences $\{\Translat\k\delta\suchthat\k\in\H/\Gamma\}$   is a
QMF-system. But its image is by what we have said again  a
QMF-system. Therefore upon setting
$$
\phi_\k = \K\,
\Translat\k\,\delta,\qquad\k\in\H/\Gamma\Topref{hhgtwzsastenm}
\Numbering
$$  we get a whole family of QMF-systems. Actually we get all of
them that way. Indeed  every unitary operator that commutes with
the
$\HH$-translations is uniquely determined by the its image of the
delta sequences $\Translat\k\delta$ with $\k\in\H/\HH$. Let us
call these functions again $\phi_\k$. We then can recover $\K$
explicitly
 by means of the following: let $\Pi_\Gamma$ be the orthogonal
projector on $\L2(\HH)$ taken in the obvious way as a subspace of
$\L2(\H)$. We now set
$$\eqalign{
\K : \s \mapsto &\sum_{\k\in\HH}
(\Translat\k\Pi_\HH\Translat\k^\ast\,\s)\,\ast\,(\Translat\k^\ast\phi_\k)
\cr &= \sum_{\k\in\HH}
(\Pi_{\Translat\k\HH}\,\s)\,\ast\,(\Translat\k^\ast\phi_\k).}
\Topref{hjkldjh}\Numbering
$$ This operator clearly commutes with the $\HH$-translations and
its image of the $\delta_\k$ with $\k\in\H/\HH$ are just the
$\phi_\k$. But now
 the $\{\phi_\k\}$ are an QMF-system and the operator as defined by
\hjkldjh\ is therefore unitary.  We therefore have shown

\Theorem{There is a bijection between QMF-systems and the group of
unitary operators $\U$ acting in $\L2(\H)$ such that
$$ [\U,\Translat\m]=0,\quad\hbox{for all $\m\in\HH$}.
$$ This bijection is explicitly given by \hhgtwzsastenm\ and
\hjkldjh.}

We now come to the Fourier space picture.    The  Fourier
transform of a sequence $\s$ over $\H=\Z^\n$ is defined as
$$
\fou\s(\omega) = \sum_{\m\in\H}
\e^{-\i\scalarproduct{\omega}{\m}},\quad
\omega\in\Torus^\n,
$$ Now consider a sub-lattice $\Gamma\subset\H$. Recall that the
dual lattice
$\Gamma^\perp$ of $\Gamma$ is defined as
$$
\k\in\Gamma^\perp\ \iff \ \e^{\i\scalarproduct{\k}{\m}}=1  \quad
\hbox{for all $\m\in\Gamma$},
$$ where the scalar-product is the one naturally inherited from
$\R^\n$.  The quotient group  $\Gamma^\perp/\H^\perp$ may be
identified with a subgroup of $\Torus^\n$.

We now want to rephrase the QMF-condition in Fourier space.
Recall that the Poisson summation formula states
 that the projection operator $\Pi_\HH$ reads in Fourier space
$$
\Pi_\HH : \fou\s(\omega) \mapsto {1\over\abs{\H/\Gamma}}
\sum_{\xi\in\HH^\perp/\H^\perp} \fou\s(\omega + \xi).
\Topref{prcxvcjkhsuoisndjf}
\Numbering
$$ As before, the  orthogonality relations can also be written
$$
\Pi_\HH\,(\tilde\phi_\k\ast\phi_\kk) =
\delta(\k-\kk)\,\delta,\qquad
\hbox{with}\ \tilde\phi_\k(\p) =\conj\phi_\k(-\p)
$$ Therefore from \prcxvcjkhsuoisndjf\ we have
$$
\sum_{\xi\in\HH^\perp/\H^\perp} \conj{\fou\phi}_\k(\omega + \xi)\,
\fou\phi_\k(\omega + \xi) = \abs{\H/\Gamma} \ \delta(\k-\kk).
$$ To put it still differently, the ortho-normality condition reads
$$
\MM_{\k,\xi}(\omega) = \inv{\sqrt{\abs{\H/\Gamma}}}\
\fou\phi_\k(\omega + \xi),
\quad\k\in\Gamma/\H, \ \xi\in\Gamma^\perp/\H^\perp
\Topref{beforegisdshnans}\Numbering
$$ is unitary for every $\omega\in\Torus^\n$.   Again the converse
is true also. To see this consider again
 the unitary operator $\K$ as given by \hjkldjh.  In Fourier space
it acts as follows
$$
\K : \fou\s(\omega) \mapsto {1\over\abs{\H/\Gamma}}\,
\sum_{\k\in\HH}
\sum_{\xi\in\HH^\perp/\H^\perp}\
\fou\s(\omega+\xi)\ \e^{i\scalarproduct{\omega+\xi}{\k}}\
\fou\phi_\k(\omega)
$$ But the functions $\fou\s(\omega+\rho)$,
$\rho\in\Gamma^\perp/\H^\perp$
 have a similar transformation behavior, namely
$$
\eqalign{
\K : \fou\s(\omega+\rho)&\mapsto{1\over\abs{\H/\Gamma}}
\sum_{\k\in\HH}
\sum_{\xi\in\Gamma^\perp/\H^\perp}\
\fou\s(\omega+\xi+\rho)\e^{i\scalarproduct{\omega+\xi+\rho}{\k}}\
\fou\phi_\k(\omega+\rho)\cr
&={1\over\abs{\H/\Gamma}}\sum_{\k\in\HH}
\sum_{\xi\in\Gamma^\perp/\H^\perp}\
\fou\s(\omega+\xi)\e^{i\scalarproduct{\omega+\xi}{\k}}\
\fou\phi_\k(\omega+\rho). }
$$ If we introduce now the vector $\Ndim$-dimensional vector with
components $\fou\s_\xi(\omega) =
\fou\s(\omega+\xi)$, $\xi\in\HH^\perp/\H^\perp$ then the action of
$\K$ becomes just matrix multiplication:
$$
\K : \s_\rho(\omega) \mapsto\sum_{\xi\in\HH^\perp/\H^\perp}
\KK_{\rho,\xi}(\omega)\s_\xi(\omega)
$$ with
$$
\KK_{\rho,\xi}(\omega)  ={1\over\abs{\H/\Gamma}} \sum_{\k\in\HH}\
\e^{i\scalarproduct{\omega+\xi}{\k}}\
\fou\phi_\k(\omega+\rho).
$$ The matrix can still be simplified by writing
$$
\KK_{\rho,\xi}(\omega) = \sum_{\k\in\HH}
\, \F_{\k,\xi} (\omega)\, \MM_{\k,\rho}(\omega)
$$ with $\MM_{\k,\rho}(\omega)$ the unitary matrix given by
\beforegisdshnans\ and
$$
\F_{\k,\xi} (\omega) = \abs{\H/\Gamma}^{-1/2}
\e^{\i\scalarproduct{\omega+\xi}{\k}}\Topref{fdfidnsdisn}\Numbering
$$
 is the complex conjugated of the QMF-matrix associated to the
QMF-system
 $\delta_\k$, $\k\in\H/\Gamma$.  Therefore in particular this
matrix is unitary  and therefore  the matrix
$\KK_{\rho,\xi}(\omega)$ is unitary.

To better understand the structure of this kind of matrix valued
function as  in \beforegisdshnans\ we introduce the permutation
representation of
$\HH^\perp/\H^\perp$; that is
$$
\Rep(\xi) = \Rep_{\omega,\rho}(\xi) = \delta_{\omega, \xi+\rho}
$$    Now the loops of the kind  \beforegisdshnans\ are precisely
the loops  $\MM(\omega) = \MM_{\k,\xi}(\omega)$ that satisfy at
$$
\MM(\omega + \rho) = \MM(\omega)\, \Rep^\ast(\rho),
\quad\rho\in\HH^\perp/\H^\perp.\Topref{thematrixeahjdiun}\Numbering
$$ The matrix valued function $\K$ defining $\U$ instead satisfies
at
$$
\K(\omega+\rho)=\Rep(\rho)\, \K(\omega)\, \Rep^\ast(\rho).
\Topref{twistingconds}
\Numbering
$$ Vice versa let  $\chi$ be a function from $\Torus^\n$ with
values in $\U(\Ndim)$
 that satisfies the above condition. Consider the associated
operator defined through its action on  the vector
$\fou\s_\xi(\omega) = \fou\s(\omega+\xi)$ for each $\omega$ via
$$
\U : \fou{\s}_\xi(\omega) \mapsto
\sum_{\rho\in\H^\perp}\chi_{\xi,\rho}(\omega)\,\fou\s_\rho(\omega)
$$ Since by hypothesis we  have $\chi(\omega+\rho) = \Rep(\rho)\,
\chi(\omega)\,\Rep^\ast(\rho)$, and hence the operation is well
defined as
 operator on a function.  In addition the  operator induced by
$\chi(\omega)$ via  its action on the $\fou\s(\omega+\xi)$  is
clearly unitary since it is invertible and it preserves the
energy.   Again the translation operator is diagonal
$$
\Translat\k : \fou\s_\xi(\omega) \to
\e^{-\i\scalarproduct{\omega}{\k}}\,
\fou\s_\xi(\omega), \quad\k\in\Gamma
$$ and thus the operator $\U$ commutes with the
$\Gamma$-translates.

\Section{Loop groups and QMF} A map $\varphi$ from the torus
$\Torus^\n$  into the group $\U(\Ndim)$ is called a loop.  The set
of loops is again a group under pointwise multiplication. All
matrix valued functions we have encountered so far are elements of
such a loop group.  The set of  loops sqtisfying at
\twistingconds\ is again a group. It  is called the group of
twisted loops. Upon conjugating a twisted loop with $\F$ defined
in \fdfidnsdisn\  we see that the twisted loops are isomorph to
the set of
$\HH^\perp$-periodic loops
$$
\eta(\omega+\rho) = \eta(\omega), \quad\rho\in\HH^\perp.
$$  This group is denoted by $\Fix(\Translat{\Gamma^\perp},
\Map(\Torus^\n\to\U(\Ndim)))$. It is again isomorphic to the whole
loop-group. Indeed by hypothesis $\HH$ is isomorphic to $\H$. Thus
$\HH^\perp$ is isomorphic to
$\H^\perp$. Let $\Phi:\R^\n\to\R^\n$ be a linear, invertible map
who sends
$\HH^\perp$  onto $\H^\perp$.  Then any loop in
$\Fix(\Translat{\Gamma^\perp},
\Map(\Torus^\n\to\U(\Ndim)))$ is mapped via
$\eta(\omega)\to\eta(\Phi(\omega))$ into a loop in
$\Map(\Torus^\n\to\U(\Ndim)))$, and this defines the isomorphism
we where looking for. As in the simplest case that we have
discussed in great detail, the QMF-matrices are right cosets of
$\Fix(\Translat{\Gamma^\perp},\Map(\Torus^\n\to\U(\Ndim)))$ and of
its analogue sub-groups respectively. The argumentation follows
exactly the same line as before and the details are left to the
reader. To sumarize we have shown

\Topref{vbcgdhjkuyfg}\Theorem{There is an explicit bijection
between
$$\eqalign{i&) \ \hbox{the QMF-filters}\cr ii&)\ \hbox{the
QMF-loops:}\
\MM(\omega\rho) = \MM(\omega)\,
\Rep^\ast(\rho),\quad\rho\in\HH^\perp\cr iii&)\ \hbox{the group of
unitary operators  with $[\K,\Shift^\H_\x]=0$,
$\x\in\HH$}\cr iv&)\ \hbox{the twisted loops}: \B(\omega\rho) =
\Rep(\rho)\,\B(\omega)\,
\Rep^\ast(\rho),\quad\rho\in\HH^\perp\cr iiv&)\ \hbox{the
untwisted loops}: \A(\omega\rho) = \A(\omega),
\quad\rho\in\HH^\perp\cr iiiv&)\ \hbox{the whole loop group} }
$$ More precisely the groups $iii$, $iv$, $iiv$ and $iiiv$ are
isomorphic.
 The set of QMF-loops is a left coset with respect to $iv$ and a
right coset with respect to $iiv$.}

We only have a closer look at the polynomial loops or what is the
same, the space of finite impulse response QMF-systems or
equivalently the subgroup  of unitary operators in $\L2(\H)$ that
commute with the $\HH$-translates and that leave invariant the
space of compactly supported sequences.

\Theorem{The group
$\Fix(\Translat{\Gamma^\perp},
\Map(\Torus^\n\to\U(\Ndim)\suchthat\hbox{poly})$ is isomorph to the
group $\Map(\Torus^\n\to\U(\Ndim)\suchthat\hbox{poly})$.}

\Proof Let $\eta\in\Fix(\Translat{\Gamma^\perp},
\Map(\Torus^\n\to\U(\Ndim))\suchthat\hbox{poly})$. Then
$$
\eta(\omega) =
\sum_{\k\in\Z}\U_\k\,\e^{\i\scalarproduct{\k}{\omega}}
$$ where the coefficients are matrices in $\U(\Ndim)$ only a
finite number of  which is different from $0$ so that the sum is
finite. In addition we have
$\eta(\omega+\xi)=\eta(\omega)$ for all $\xi\in\HH^\perp$. This
implies
$$
\left\{\sum_{\k\in\HH} + \sum_{\k\not\in\HH}\right\}
\U_\k\,\e^{\i\scalarproduct{\k}{\omega + \xi}} =
\left\{\sum_{\k\in\HH} + \sum_{\k\not\in\HH}\right\}
 \U_\k\,\e^{\i\scalarproduct{\k}{\omega}}.
$$ Since $\e^{\i\scalarproduct{\k}{\xi}} =1$ in the first sum, it
follows that
$$
\sum_{\k\not\in\HH}\, \U_\k\,(\e^{\i\scalarproduct{\k}{\xi}}-1)\,
\e^{\i\scalarproduct{\k}{\omega}} = 0,
$$ which implies, since $\e^{\i\scalarproduct{\k}{\xi}}\not=1$ in
this sum, that
$\U_\k=0$ for $\k\not\in\HH$. Therefore the group isomorphism
defined above  gives the desired isomorphism for polynomial loops.
\qed

\Subsection{The factorization problem.} In more than one dimension
it really is a problem. Only in the case $\n=1$ and arbitrary
$\Ndim$ factorization theorems for the polynomial loops are known.
In all other cases the problem is open. However there is a general
construction to obtain a large quantity  of polynomial loops. It
uses once more the equivalence between the  unitary operators
acting in $\L2(\H)$ and the loop groups. As in the easiest case we
may define a unitary operator $\U_\u^\F$ whose restriction to some
fundamental domain $\F=\HH/\H$---consisting of $\Ndim$
non-congruent points modulo $\HH$---coincides with the unitary map
$\u\in\U(\Ndim)$. This defines a whole manifold of finite impulse
response QMF-systems and hence of polynomial loops. By changing
the fundamental domain we obtain another family of such unitary
operators. In general the operators corresponding to different
fundamental domains do not commute and thus one can compose them
to obtain nontrivial new ones. The kind of factorization we
propose therefore is the same as already proposed in
\def\hols{[3]}\hols.
$$
\U = \prod_{\p=0}^\m\,\U_{\u_\p}^{\F_\p},\quad
\hbox{with}\quad\u_\p\in\U(\Ndim),\ \F_\p\ \hbox{some fundamental
domaine}
$$ whether or not this family is exhaustive is not known yet. In
any case it is a subgroup of all polynomial  loops.

\Section{References}
\Ref{Estaban D. Galand C. (1977), {\it Application of quadrature
mirror
 filters to split-band voice coding schemes},  Proc. IEEE Int.
Conf. Acoust. Signal Speech Process., Hartford, Connecticut, pp.
191-195}

\Ref{Daubechies I. {\it Ten Lectures on Wavelets}, CBMS-NSF
Lecture Notes nr. 61, SIAM, Philadelphia (1992)}

\Ref{Holschneider M. (1992),  {\it Wavelet analysis over abelian
groups},
 ACHA, to appear}

\Ref{Pollen D., {\it $SU_I(2,F[z,1/z])$ for $F$ a subfield of
$\CC.$},  JAMS, 3, 1990. }

\Ref{Pressley A. \&\ Segal G.  (1986), {\it Loop Groups}, Oxford
Mathematical
 Monograph}

\bye